%% file: main.tex
\documentclass[11pt, reqno, psamsfonts]{amsart}
\usepackage{format}

\input{author}
\title{\sffamily Borel Local Lemma: arbitrary random variables and \\ limited exponential growth}
\begin{document}

    \vspace*{-0pt}

    \maketitle    

    
    \begin{abstract}
        The Lov\'asz Local Lemma (the LLL for short) is a powerful tool in probabilistic combinatorics that is used to verify the existence of combinatorial objects with desirable properties. Recent years saw the development of various ``constructive'' versions of the LLL. A major success of this research direction is the Borel version of the LLL due to  Cs\'oka, Grabowski, M\'ath\'e, Pikhurko, and Tyros, which holds under a subexponential growth assumption. 
        A drawback of their approach is that it only applies when the underlying random variables take values in a finite set. We present an alternative proof of a Borel version of the LLL that holds even if the underlying random variables are continuous and applies to dependency graphs of limited exponential growth.
    \end{abstract}

\input{introduction}

\input{MT}

\subsection*{Acknowledgments}

    We thank Tam\'as K\'atay, Andrew Marks, and the anonymous referees for their helpful comments. This material is based upon work partially supported by the Alfred P. Sloan Foundation, the National Science Foundation under grant DMS-2528522, and the National Natural Science Foundation of China grants 12371343, 12525110. Any opinions, findings, and conclusions or recommendations expressed in this material are those of the authors and do not necessarily reflect the views of the funding agencies.

\printbibliography
\end{document}

%% file: author.tex
\author{Anton~Bernshteyn}
\address{\normalfont (AB) Department of Mathematics, University of California, Los Angeles, CA, USA}
\email{bernshteyn@math.ucla.edu}

\author{Jing~Yu}
\address{\normalfont (JY) Shanghai Center for Mathematical Sciences, Fudan University, Shanghai, China}
\email{jyu@fudan.edu.cn}

\thanks{AB's research is partially supported by the NSF CAREER grant DMS-2239187 and the Sloan Research Fellowship (2025). JY's research is partially supported by the National Natural Science Foundation of China grants 12371343 and 12525110 (PI: Hehui Wu).}

%% file: introduction.tex
\section{Introduction}

    This paper is a contribution to \emph{Borel combinatorics}: an area at the crossroads of combinatorics and descriptive set theory that aims to perform \emph{combinatorial} constructions using \emph{Borel} sets and functions. For an overview of this subject, see the survey \cite{KechrisMarks} by Kechris and Marks and introductory articles \cite{pikhurko2020borel} by Pikhurko and \cite{bernshteyn2022descriptive} by the first named author. We also direct the reader to \cite{KechrisDST,AnushDST} for a general introduction to descriptive set theory. While most work in Borel combinatorics addresses specific combinatorial problems, a recent trend is to design general tools that can be applied to broad problem classes \cite{bernshteyn2019measurable, Ber_Baire, BernshteynDistributed, bernshteyn2023probabilistic, trees, CGMPT, paths, GRgrids, ASIalgorithms, FelixTrees,BerFelixASI}. One such tool is the \emph{Lov\'asz Local Lemma} (the \emph{LLL} for short): a probabilistic statement that, under certain numerical conditions, guarantees the existence of a labeling satisfying a set of constraints \cites{EL}{SpencerLLL}[Chap.~5]{AS}. Recent years saw an explosion of ``constructive'' versions of the LLL, for various notions of ``constructiveness'' \cite{Beck,BGRDeterministicLLL,BMUDeterministicLLL,FG,MT,RSh,bernshteyn2019measurable,BernshteynDistributed, bernshteyn2023probabilistic, CGMPT,BerFelixASI}. Here we are interested in the behavior of the LLL in the Borel setting. Specifically, our goal is to generalize and sharpen the Borel LLL for dependency graphs of subexponential growth due to Cs\'oka, Grabowski, M\'ath\'e, Pikhurko, and Tyros \cite{CGMPT}.

    There exist several essentially equivalent frameworks for stating the LLL in the Borel context \cite{bernshteyn2019measurable,bernshteyn2022descriptive,CGMPT}. We shall find the following formalism convenient:

    \begin{defn}[Constraint satisfaction problems]\label{defn:CSP}
        A \emphd{constraint satisfaction problem} (a \emphd{CSP} for short) is a tuple $\Pi = (V, \Lambda, \Const, \dom, \Bad)$, where:
        \begin{itemize}
            \item $V$, $\Lambda$, and $\Const$ are sets whose elements are called \emphd{variables}, \emphd{labels}, and \emphd{constraints} respectively,

            \item 
            $\dom$ assigns to each constraint $\const \in \Const$ a finite set $\dom(\const) \subseteq V$ called the \emphd{domain} of $\const$,

            \item $\Bad$ assigns to each constraint $\const \in \Const$ a set $\Bad(\const) \subseteq \Lambda^{\dom(\const)}$ of \emphd{bad labelings} $\phi \colon \dom(\const) \to \Lambda$.
        \end{itemize}
        A function $f \colon V \to \Lambda$ is called a \emphd{labeling} of $V$. Such a labeling $f$ \emphd{violates} a constraint $\const \in \Const$ if its restriction to $\dom(\const)$ is in $\Bad(\const)$; otherwise, $f$ \emphd{satisfies} $\const$.  A \emphd{solution} to $\Pi$ is a labeling $f \colon V \to \Lambda$ that satisfies every constraint $\const \in \Const$. 
    \end{defn}

   Let us illustrate this definition with a couple of simple examples showing how natural combinatorial problems can be represented by CSPs.

    \begin{exmp}[Proper coloring]
        Given a graph $G$ and a natural number\footnote{In this paper, $0$ is a natural number, and a natural number $k$ is identified with the $k$-element set $\set{i \in \N \,:\, i < k}$.} $k \in \N$, a \emphd{proper $k$-coloring} of $G$ is a function $f \colon V(G) \to k$ such that $f(u) \neq f(v)$ whenever $\set{u,v} \in E(G)$. It is clear that a proper $k$-coloring of $G$ is the same as a solution to the CSP
        \[
            \Pi_k(G) \,\defeq\, \big(V(G), \, k,\, E(G), \, \mathsf{id} , \,  \Bad\big),
        \]
        where 
        for each edge $e \in E(G)$, $\mathsf{id}(e) \defeq e$ (i.e., the domain of the corresponding constraint is 
        the set containing the two endpoints of $e$), while $\Bad(e)$ is the set of all constant mappings $\phi \colon e \to k$.
    \end{exmp}

    \begin{exmp}[Sinkless orientation]\label{exmp:sinkless}
        Let $G$ be a locally finite graph. An orientation of $G$ is \emphd{sinkless} if it has no \emphd{sinks}, i.e., no vertices with outdegree $0$. The problem of finding a sinkless orientation of $G$ can be encoded by a CSP $\Pi_{\bullet\to}(G)$ as follows. Fix an arbitrary orientation $\vec{G}_0$ of $G$. We can then identify any orientation $\vec{G}$ of $G$ with a labeling $f \colon E(G) \to \set{{+},{-}}$, where $f(e) = {+}$ if and only if $e$ is oriented the same in $\vec{G}$ and $\vec{G}_0$. Now we let
        \[
            \Pi_{\bullet\to}(G) \,\defeq\, \big(E(G), \, \set{{+},{-}}, \, V(G), \, \partial_G, \, \Bad\big),
        \]
        %
        where for each vertex $v \in V(G)$, $
        \partial_G(v)$ 
        is the set of all edges of $G$ incident to $v$ and $\Bad(v)$ contains the unique mapping $\phi \colon \partial_G(v) \to \set{{+},{-}}$ that makes $v$ a sink.
    \end{exmp}

    The list of examples can be continued almost indefinitely and includes such classical problems as perfect matching, maximal independent set, edge coloring, hypergraph coloring, graph homomorphism, and so on. We remark that in most combinatorial applications, the set of labels $\Lambda$ is finite; however, that is not a requirement in Definition~\ref{defn:CSP}.

    In view of the versatility of Definition~\ref{defn:CSP}, it is desirable to develop flexible general conditions that guarantee the existence of a solution to a given CSP $\Pi$. 
    The LLL gives just such a condition. Roughly, the LLL says that $\Pi$ has a solution provided that:
    \begin{enumerate*}[label=\ep{\itshape\alph*}]
        \item the constraints interact with each other ``in a limited way,'' and
        \item each individual constraint is ``unlikely'' to be violated.
    \end{enumerate*}
    
    The 
    interaction between the constraints is captured by 
    the \emph{dependency graph} of $\Pi$:

    \begin{defn}[Dependency graph]
        Let $\Pi = (V, \Lambda, \Const, \dom, \Bad)$ be a CSP. The \emphd{dependency graph} of $\Pi$, denoted by $D_\Pi$, is the graph with vertex set $\Const$ and edge set
        $
            \set{\set{\const,\const'} \,:\, \const \neq \const', \,  \dom(\const) \cap \dom(\const') \neq \0}.
        $
    \end{defn}

    Throughout, we shall assume the dependency graph $D_\Pi$ has finite maximum degree. Intuitively, the smaller the maximum degree of $D_\Pi$ is, the less the constraints interact with each other. 

    To measure how likely each constraint is to be violated, we equip the label set $\Lambda$ with a probability measure $\P$. To unclutter the notation, we shall, when there is no possibility of confusion, use the symbol ``$\P$'' not only for the given measure on $\Lambda$ but also for product measures of the form $\P^S$ on $\Lambda^S$ for some set $S$. Now, for a constraint $\const \in \Const$, we can randomly generate a labeling $\phi \colon \dom(\const) \to \Lambda$ (i.e., a point in $\Lambda^{\dom(\const)}$) by drawing the labels $\phi(v)$ for each $v \in \dom(\const)$ independently from $(\Lambda, \P)$. The likelihood of violating $\const$ is simply the probability the resulting labeling is bad, i.e., $\P[\Bad(\const)]$. 


    With these preliminaries, we are ready to state the Lov\'asz Local Lemma:

    \begin{theo}[{Lov\'asz Local Lemma \cites{EL}{SpencerLLL}[Corollary 5.1.2]{AS}}]\label{theo:LLL}
        Let $\Pi= (V, \Lambda, \Const, \dom, \Bad)$ be a CSP and let $\P$ be a probability measure on $\Lambda$ making $(\Lambda, \P)$ a standard probability space. 
        If there exist $p \in [0,1)$ and $d \in \N$ such that:
        \begin{itemize}
            \item $\Bad(\const)$ is $\P$-measurable and $\P[\Bad(\const)] \leq p$ for all $\const \in \Const$,
            \item the maximum degree of the dependency graph $D_\Pi$ is at most $d$, and
            \item $\e  \, p\, (d+1) < 1$ \ep{where $\e  = 2.71\ldots$ is the base of the natural logarithm},
        \end{itemize}
        then $\Pi$ has a solution $f \colon V \to \Lambda$.
    \end{theo}

    \begin{remks}\label{remk:LLL}
        \begin{enumerate}[wide, label=\ep{\normalfont\roman*}]
        \item The above statement of the LLL is somewhat restricted compared to the more standard version typically found in the combinatorics literature, in that we compute the probability $\P[\Bad(\const)]$ with respect to specifically the \emph{product} measure on $\Lambda^{\dom(\const)}$. 
        This 
        set-up is sometimes called the \emph{variable version} of the LLL. Even though this setting is not the most general, it does encompass virtually all standard applications \cite[41]{MolloyReed} and is often viewed as the ``right one'' for algorithmic considerations (see, e.g., \cite{Beck,MT,FG}). That being said, it is possible to formulate the LLL without assuming a product structure on the underlying probability space; see \cite[\S5.1]{AS} for the statement of the LLL in abstract probability spaces and \cite{KolipakaSzegedy,Experimental,VLLL} for comparisons between the abstract and variable versions of the LLL.

        \item\label{item:compactness} Although the LLL is usually stated and proved in the case when the set of constraints $\Const$ is finite, Theorem~\ref{theo:LLL} remains valid even if $\Const$ is infinite. Indeed, we may without loss of generality assume that $\Lambda$ is a compact metric space and $\P$ is a Borel probability measure \cite[Thm.~17.41]{KechrisDST}. Pick any $p' > p$ such that $\e \,p'\, (d+1) < 1$ and replace each $\Bad(\const)$ by an \emph{open} set $\Bad'(\const) \supseteq \Bad(\const)$ of measure at most $p'$, which is possible since $\P$ is regular \cite[Thm.~17.10]{KechrisDST}. The existence of a solution to $\Pi$ now follows by applying the LLL to finite subsets of $\Const$ and using the compactness of the space $\Lambda^V$; see, e.g., \cite[proof of Thm.~5.2.2]{AS} for details.
    \end{enumerate}
    \end{remks}

    In this paper we are interested in the following general question:

    \begin{ques}
        Let $\Pi = (V, \Lambda, \Const, \dom, \Bad)$ be a CSP such that $V$ and $\Lambda$ are standard Borel spaces. When can the LLL be invoked to conclude that there exists a \emph{Borel} solution $f \colon V \to \Lambda$ to $\Pi$?
    \end{ques}

    Recall that a \emphd{standard Borel space} is a set $X$ equipped with a $\sigma$-algebra $\mathfrak{B}(X)$ generated by a complete separable metric on $X$. The members of $\mathfrak{B}(X)$ are called the \emphd{Borel subsets} of $X$. By the Borel Isomorphism Theorem \cite[Thm.~15.6]{KechrisDST}, all uncountable standard Borel spaces are isomorphic to each other, so, as a rule, no generality is lost by replacing $X$ with some familiar space, such as the real line $\R$. A function $f \colon X \to Y$ between standard Borel spaces is \emphd{Borel} if $f$-preimages of Borel subsets of $Y$ are Borel in $X$, or, equivalently, if the graph of $f$ is a Borel subset of $X \times Y$ \cite[Thm.~14.12]{KechrisDST}. Borel sets and functions have a variety of useful regularity properties and are generally taken to represent the ``gold standard'' of well-behavedness in descriptive set theory.

    Unfortunately, as the next example illustrates, even under very favorable circumstances, the solutions produced by the LLL may fail to be Borel.

    \begin{exmp}[Borel sinkless orientation---a negative result]\label{exmp:sinkless_negative}
        Fix $d \in \N$ and let $G$ be a $d$-regular graph \ep{meaning that every vertex of $G$ has exactly $d$ neighbors}. Recall the sinkless orientation problem $\Pi_{\bullet\to}(G)$ from Example~\ref{exmp:sinkless}. It is easy to see that the dependency graph of $\Pi_{\bullet \to}(G)$ is $G$ itself, hence its maximum degree is $d$. Letting $\P$ be the uniform probability distribution on $\set{+,-}$, we have $\P[\Bad(v)] = 2^{-d} \eqqcolon p$ for each $v \in V(G)$, because exactly one out of the $2^d$ ways to orient the edges incident to $v$ is bad. The inequality
        \[
            \e\, p \, (d+1) \,=\, \e \, 2^{-d} \, (d+1) \,<\, 1
        \]
        is satisfied for all $d \geq 4$; moreover, its left-hand side rapidly approaches $0$ as $d$ goes to infinity. Nevertheless, no matter how large $d$ is, $G$ may fail to have a Borel sinkless orientation, since Thornton \cite[Thm.~3.5]{thornton2022orienting}, using the determinacy method of Marks \cite{marks2016determinacy}, constructed for each $d \in \N$ a Borel $d$-regular graph $G$ with no Borel sinkless orientation.
    \end{exmp}

    While the above example shows the LLL may fail to produce a Borel solution in general, there has been a recent string of positive ``LLL-like'' results in descriptive set theory \cite{bernshteyn2019measurable,bernshteyn2022descriptive,CGMPT,BerFelixASI,bernshteyn2023probabilistic}. Most of these results yield weaker conclusions (e.g., measurable in place of Borel) or require stronger bounds than $\e \, p\, (d+1) < 1$ (often both). A notable exception is the Borel version of the LLL due to Cs\'oka, Grabowski, M\'ath\'e, Pikhurko, and Tyros \cite{CGMPT}. Their result yields a Borel solution to a given CSP $\Pi$ under the usual LLL condition $\e  \, p \, (d+1) < 1$; to achieve this, they add some extra assumptions on $\Pi$, the main one being that the dependency graph $D_\Pi$ ought to be \emph{of subexponential growth}.

    \begin{defn}[Growth of graphs]
        The \emphd{growth function} $\gamma_G \colon \N \to \N \cup \set{\infty}$ of a locally finite graph $G$ is 
        given by $\gamma_G(R) \defeq \sup_{v \in V(G)} |B_G(v, R)|$, where $B_G(v,R)$ is the \emphd{$R$-ball around $v$ in $G$}, i.e., the set of all vertices joined to $v$ by a path of at most $R$ edges. 
        The quantity 
        \[
            \egr(G) \,\defeq\, \lim_{R \to \infty} \sqrt[R]{\gamma_G(R)} \,=\, \inf_{R \geq 1} \sqrt[R]{\gamma_G(R)}
        \]
        is called the \emphd{exponential growth rate} of $G$. 
        (The limit exists and is equal to the infimum by Fekete's lemma \cite[Lem.~A.4.2]{FracBook}.) 
        We say that $G$ is \emphd{of subexponential growth} if $\egr(G) = 1$.
    \end{defn}

    The theorem of Cs\'oka \emph{et al.}~and all the results of this paper apply to \emph{Borel CSPs}, i.e., CSPs of the form $\Pi = (V, \Lambda, \Const, \dom, \Bad)$ in which $V$, $\Lambda$, and $\Const$ are standard Borel spaces and the assignments $\dom$ and $\Bad$ are Borel in a certain natural sense; see \S\ref{subsec:Borel_defns} for the formal definition. Cs\'oka \emph{et al.}~showed that a Borel CSP $\Pi$ that fulfills the assumptions of the LLL has a Borel solution, provided that its dependency graph is of subexponential growth and, furthermore,  the label set $\Lambda$  and
         the \emph{order} of $\Pi$, defined below, are finite. 


    \begin{defn}[Order and maximum variable degree]
        Let $\Pi=(V, \Lambda, \Const, \dom, \Bad)$ be a CSP. For $v \in V$, let $\dom^{-1}(v) \defeq \set{\const \in \Const \,:\, v \in \dom(\const)}$.  The \emphd{order} and the \emphd{maximum variable degree} of $\Pi$ are 
        \[
            \ord(\Pi) \,\defeq\, \sup_{\const \in \Const}|\dom(\const)| \qquad \text{and} \qquad \vdeg(\Pi) \,\defeq\, \sup_{v  \in V} |\dom^{-1}(v)|. 
        \]
    \end{defn}

    \begin{theo}[{Cs\'oka--Grabowski--M\'ath\'e--Pikhurko--Tyros \cite{CGMPT}}]\label{theo:CGMPT}
        Let $\Pi = (V, \Lambda, \Const, \dom, \Bad)$ be a Borel CSP and let $\P$ be a Borel probability measure on $\Lambda$. If there are $p \in [0,1)$, $d \in \N$ such that:
        \begin{enumerate}[label={\normalfont\ep{\roman*}}]
            \item\label{item:LLL} $\P[\Bad(\const)] \leq p$ for all $\const \in \Const$, the maximum degree of $D_\Pi$ is at most $d$, and $\e \, p\, (d+1) < 1$,
                \item\label{item:subexp} $D_\Pi$ is of subexponential growth, and
            \item\label{item:finite} $\Lambda$ is a finite set and $\ord(\Pi) < \infty$,
        \end{enumerate}
        then $\Pi$ has a Borel solution $f \colon V \to \Lambda$.
    \end{theo}

    In Theorem~\ref{theo:CGMPT}, \ref{item:LLL} is just the standard LLL assumption (as in Theorem~\ref{theo:LLL}). Example~\ref{exmp:sinkless} shows that some extra conditions are necessary to guarantee the existence of a Borel solution, and the subexponential growth of $D_\Pi$ is a natural requirement to add. Indeed, Theorem~\ref{theo:CGMPT} belongs to a growing body of research showing that various combinatorial problems can be solved in a Borel way on graphs of subexponential growth \cite{ConleyTamuz, thornton2022orienting,SubexpVizing}. For example, Thornton proved that if $d \geq 3$, then every $d$-regular Borel graph of subexponential growth has a Borel sinkless orientation \cite[Thm.~1.5]{thornton2022orienting}. For $d \geq 4$, this fact is a consequence of Theorem~\ref{theo:CGMPT}. Other applications of Theorem~\ref{theo:CGMPT} can be found in \cite{BY23, BernshteynDistributed}.

    On the other hand, condition~\ref{item:finite} in Theorem~\ref{theo:CGMPT} may seem somewhat puzzling. While it is usually satisfied in combinatorial applications, it is unclear whether its presence is necessary, given that the classical LLL (i.e., Theorem~\ref{theo:LLL}) is valid even if $(\Lambda, \P)$ is a \emph{continuous} probability space. As we explain below, the proof of Theorem~\ref{theo:CGMPT} in \cite{CGMPT} uses the finiteness of both $\Lambda$ and $\ord(\Pi)$ in an essential and apparently unavoidable way. Nevertheless, in this paper we present an alternative proof strategy that completely eliminates assumption \ref{item:finite}:

    \begin{tcolorbox}
    \begin{theo}\label{theo:main_subexp}
        Let $\Pi = (V, \Lambda, \Const, \dom, \Bad)$ be a Borel CSP and let $\P$ be a Borel probability measure on $\Lambda$. Suppose there exist $p \in [0,1)$ and $d \in \N$ such that:
        \begin{enumerate}[label={\normalfont\ep{\roman*}}]
            \item $\P[\Bad(\const)] \leq p$ for all $\const \in \Const$, the maximum degree of $D_\Pi$ is at most $d$, and $\e \, p\, (d+1) < 1$,
                \item $D_\Pi$ is of subexponential growth.
        \end{enumerate}
        Then $\Pi$ has a Borel solution $f \colon V \to \Lambda$.
    \end{theo}
    \end{tcolorbox}


    Theorem~\ref{theo:CGMPT} is a corollary of the following more general statement established in \cite{CGMPT}, which explicitly invokes the finite values of $|\Lambda|$ and $\ord(\Pi)$: 

    \begin{theo}[{\cite[Thm.~4.5]{CGMPT}}]\label{theo:CGMPT_slow}
        Let $\Pi = (V, \Lambda, \Const, \dom, \Bad)$ be a Borel CSP such that $\Lambda$ is a nonempty finite set, and let $\P$ be the uniform probability distribution on $\Lambda$. 
        Suppose there exist parameters $p \in [0,1)$, $d$, $\Delta \in \N$, and $\epsilon > 0$ such that:
        \begin{enumerate}[label=\ep{\normalfont\roman*}]
            \item $\P[\Bad(\const)] \leq p$ for all $\const \in \Const$ and the maximum degree of $D_\Pi$ is at most $d$,
            \item\label{item:Delta} both $\ord(\Pi)$ and $\vdeg(\Pi)$ are at most $\Delta$, 
            \item\label{item:growth_rate} $\egr(D_\Pi) < (1+\epsilon)^{2/3}$, and 
            \item\label{item:new_bound} $\e  \, p\, (d+1) \, |\Lambda|^{\epsilon \Delta} < 1$.
        \end{enumerate}
        Then $\Pi$ has a Borel solution $f \colon V \to \Lambda$.
\end{theo}

    \begin{remks}
        \begin{enumerate}[wide, label=\ep{\normalfont\roman*}]
        \item Since our notation is somewhat different from that in \cite{CGMPT}, let us briefly comment on how Theorem~\ref{theo:CGMPT_slow} follows from \cite[Thm.~4.5]{CGMPT}. Let $G$ be the bipartite graph with a bipartition $(V, \Const)$ such that $v \in V$ is adjacent to $\const \in \Const$ if and only if $v \in \dom(\const)$ (this is the construction in \cite[Rmk.~1.3]{CGMPT}). Condition \ref{item:Delta} of Theorem~\ref{theo:CGMPT_slow} then says that the maximum degree of $G$ is at most $\Delta$. It is not hard to see that if $\ord(\Pi) < \infty$, then $\egr(D_\Pi) = \egr(G)^2$, so \ref{item:growth_rate} is equivalent to $\egr(G) < (1+\epsilon)^{1/3}$. This implies that for large enough $R \in \N$, $\gamma_{G}(3R) < (1+\epsilon)^R$. In the notation of \cite{CGMPT}, this means that $G \in \mathrm{SubExp}(R, \epsilon,\Delta)$, and hence \cite[Thm.~4.5]{CGMPT} may be applied to obtain a Borel solution to $\Pi$.

        \item Let us also sketch how Theorem~\ref{theo:CGMPT_slow} yields Theorem~\ref{theo:CGMPT}. Suppose we are in the setting of Theorem~\ref{theo:CGMPT}. To begin with, we may replace $\Pi$ by an equivalent problem (possibly with a larger set of labels) to arrange that $\P$ is the uniform probability measure on $\Lambda$ \cite[Rmk.~2.4]{BY23}. Next we note that $\Delta \defeq \max \big\{\vdeg(\Pi), \ord(\Pi)\big\}$ must be finite, because $\ord(\Pi) < \infty$ by assumption and $\vdeg(\Pi) \leq d+1$. Hence, we can choose $\epsilon > 0$ so small that $\e  \, p\, (d+1)< |\Lambda|^{-\epsilon \Delta}$. Since $D_\Pi$ is of subexponential growth, $\egr(D_\Pi) = 1 < (1+\epsilon)^{2/3}$, so $\Pi$ has a Borel solution by Theorem~\ref{theo:CGMPT_slow}.
    \end{enumerate}
    \end{remks}

    Note that Theorem~\ref{theo:CGMPT_slow} holds even if the dependency graph $D_\Pi$ is \emph{of exponential growth}, as long as its exponential growth rate is sufficiently small. Unfortunately, how small the growth rate has to be significantly depends on $|\Lambda|$ and $\Delta$; that is why both $\Lambda$ and $\ord(\Pi)$ must be finite for this result to apply. In particular, notice that Theorem~\ref{theo:CGMPT_slow} can only be used when $\egr(D_\Pi) < 2^{2/3} \approx 1.59$. Indeed, in the setting of Theorem~\ref{theo:CGMPT_slow}, every constraint $\const \in \Const$ must satisfy either $\Bad(\const) = \0$ or
    \[
        \P[\Bad(\const)] \,\geq\, \frac{1}{|\Lambda|^{|\dom(\const)|}} \,\geq\, \frac{1}{|\Lambda|^{\ord(\Pi)}} \,\geq\, \frac{1}{|\Lambda|^\Delta}.
    \]
    Therefore, except for the trivial case $p = 0$, we have $p \geq |\Lambda|^{-\Delta}$. Hence, the bound $\e  \, p\, (d+1) \, |\Lambda|^{\epsilon \Delta} < 1$ can only be satisfied if $\epsilon < 1$, and thus Theorem~\ref{theo:CGMPT_slow} requires $\egr(D_\Pi) < 2^{2/3}$.

    In our main result we remedy the above shortcomings of Theorem~\ref{theo:CGMPT_slow}; that is, we generalize it to arbitrary $\Lambda$, remove the bound on $\ord(\Pi)$, and make it applicable to dependency graphs whose exponential growth rate may be arbitrarily large (although it must still be bounded in terms of the relationship between $p$ and $d$):

    \begin{tcolorbox}
    \begin{theo}\label{theo:main}
        Let $\Pi = (V, \Lambda, \Const, \dom, \Bad)$ be a Borel CSP and let $\P$ be a Borel probability measure on $\Lambda$. Suppose there exist parameters $p \in [0,1)$, $d \in \N$, and $s > 1$ such that:
        \begin{enumerate}[label=\ep{\normalfont\roman*}]
            \item $\P[\Bad(\const)] \leq p$ for all $\const \in \Const$ and the maximum degree of $D_\Pi$ is at most $d$,
            \item $\egr(D_\Pi) < s$, and
            \item\label{item:new_bound} $p\, (\e \,(d+1))^{s} < 1$.
        \end{enumerate}
        Then $\Pi$ has a Borel solution $f \colon V \to \Lambda$.
    \end{theo}
    \end{tcolorbox}

    Theorem~\ref{theo:main_subexp} is an immediate corollary to Theorem~\ref{theo:main}: if $\egr(D_\Pi) = 1$ and $\e \, p\, (d+1) < 1$, we can find $\epsilon > 0$ such that $p\, (\e\,(d+1))^{1+\epsilon} < 1$ and apply Theorem~\ref{theo:main} with $s = 1 + \epsilon$. Furthermore, it is a standard observation (see, e.g., \cite{Beck, bernshteyn2022descriptive, CP, FG, GHK,BerFelixASI}) that 
     for a majority of applications, instead of the usual LLL condition $\e\,p\, (d+1) < 1$, it is enough to have a version of the LLL that holds under a \emph{polynomial criterion}, i.e., with a bound of the form $p \, f(d) < 1$ for some polynomial $f$. When $s$, the bound on the exponential growth rate of $D_\Pi$, is treated as a constant parameter, Theorem~\ref{theo:main} becomes of this type.

     \begin{exmp}[Borel sinkless orientation---a positive result]
         Consider again the sinkless orientation problem $\Pi_{\bullet\to}(G)$ for a $d$-regular Borel graph $G$ from  Examples~\ref{exmp:sinkless} and \ref{exmp:sinkless_negative}. It is easy to see that if every component of $G$ has a cycle, then $G$ has a Borel sinkless orientation. Furthermore, when $G$ is of subexponential growth, Thornton showed that $G$ has a Borel orientation in which every vertex has outdegree at least $d/2 - 1$ \cite[Thm.~1.5]{thornton2022orienting}. Theorem~\ref{theo:main} yields an extension of this result to graphs whose exponential growth rate is near-linear in $d$: For any $\epsilon > 0$, $G$ has a Borel orientation where every vertex has outdegree at least $(1-\epsilon)d/2$, provided that $\egr(G) \leq O(\epsilon^2 d/\log d)$.  
     \end{exmp}

    In addition to establishing a more general result, our proof of Theorem~\ref{theo:main} is also somewhat simpler than the proof of Theorem~\ref{theo:CGMPT_slow} given in \cite{CGMPT}, although the two arguments have several common ingredients. The starting point of both approaches is the so-called \emph{Moser--Tardos Algorithm}: a randomized procedure for solving CSPs developed and analyzed by Moser and Tardos in their landmark paper \cite{MT}. We give an overview of the Moser--Tardos Algorithm in \S\ref{subsec:MTA}. The difficulty in the Borel setting, roughly, is that the Moser--Tardos Algorithm needs a large set of mutually independent random inputs, which cannot be generated in a Borel way. The authors of \cite{CGMPT} deal with this challenge via \emph{randomness conservation}: namely, they reuse the same random input multiple times. With a careful analysis of the Moser--Tardos Algorithm, they are able to show that it still succeeds when randomness is conserved in this way. This step in the proof involves combinatorial counting arguments that are ultimately responsible for the dependence on $|\Lambda|$ and $\Delta$ in Theorem~\ref{theo:CGMPT_slow}. Our strategy is to use \emph{probability boosting} instead: thanks to the bound on the exponential growth rate of $D_\Pi$, we are able to reduce the given CSP $\Pi$ to a different CSP $\Pi'$ whose corresponding value of $p$ is \emph{much} smaller. This reduction is explained in \S\ref{subsec:good_definitions}, and the desired properties of $\Pi'$ are established in \S\S\ref{subsec:facts}, \ref{subsec:good}, and \ref{subsec:witness}. We then find a Borel solution to $\Pi'$---and hence to $\Pi$---using the following lemma:

    \begin{Lemma}\label{lemma:double_exp}
        Let $\Pi = (V, \Lambda, \Const, \dom, \Bad)$ be a Borel CSP and let $\P$ be a Borel probability measure on $\Lambda$. Suppose there exist parameters $p \in [0,1)$ and $d \in \N$ such that:
        \begin{itemize}
            \item $\P[\Bad(\const)] \leq p$ for all $\const \in \Const$, the maximum degree of $D_\Pi$ is at most $d$, and
            \[p\, (d+1)^{d+1} \,<\, 1.\]
        \end{itemize}
        Then $\Pi$ has a Borel solution $f \colon V \to \Lambda$.
    \end{Lemma}

    Lemma~\ref{lemma:double_exp} is proved via the \emph{method of conditional probabilities}, a standard derandomization technique in computer science \cites[\S16]{AS}[\S5.6]{RandAlg}. This method has already been applied to obtain ``constructive'' versions of the LLL in various contexts \cite{Beck,FG,BerFelixASI,bernshteyn2023probabilistic}. Indeed, the proof of Lemma~\ref{lemma:double_exp} is almost exactly the same as that of \cite[Thm.~1.6]{bernshteyn2023probabilistic} and \cite[Thm.~3.3]{BerFelixASI}. The only difference is that the results in \cite{bernshteyn2023probabilistic,BerFelixASI} are stated for finite $\Lambda$; however, their proofs can be adapted for the general case with minimal modifications. For completeness, we present a direct proof of Lemma~\ref{lemma:double_exp} in \S\ref{subsec:proof_de}. We then finish the proof of Theorem~\ref{theo:main} in \S\ref{subsec:final}.

%% file: MT.tex
\section{The Moser--Tardos Algorithm and its local analysis}\label{sec:MT}

    In this section we present the combinatorial heart of our argument. Throughout \S\ref{sec:MT}, we fix a CSP $\Pi=(V, \Lambda, \Const, \dom, \Bad)$ and a probability measure $\P$ on $\Lambda$ making $(\Lambda, \P)$ a standard probability space. We also let $D \defeq D_\Pi$ be the dependency graph of $\Pi$ and assume that $D$ is locally finite. For each $\const \in \Const$, we write $N[\const] \defeq B_D(\const, 1)$ and $N(\const) \defeq N[\const] \setminus \set{\const}$ for the \emphd{closed}, resp.~\emphd{open} \emphd{neighborhood} of $\const$ in the graph $D$. Note that since $D$ is locally finite by assumption, $N[\const]$ and $N(\const)$ are finite sets. 

\subsection{The Moser--Tardos Algorithm}\label{subsec:MTA}

    As mentioned in the introduction, an important tool in the study of ``constructive'' aspects of the LLL is the \emphd{Moser--Tardos Algorithm}, or $\MTA$ for short (see Algorithm~\ref{alg:MT}).
    It was introduced by Moser and Tardos in \cite{MT} in order to prove an algorithmic version of the LLL. $\MTA$ takes as input a CSP $\Pi$ and a map $\tbl \colon V \to \Lambda^\N$, called a \emphd{table} \ep{we picture it as a matrix with rows indexed by $\N$, columns indexed by $V$, and entries from $\Lambda$}. For convenience, given $v \in V$ and $n \in \N$, we write $\tbl(v,n)$ for the $n$-th entry of the sequence $\tbl(v) \in \Lambda^\N$. The algorithm attempts to solve $\Pi$ by building a sequence of labelings $f_0$, $f_1$, \ldots{} $\colon V \to \Lambda$ as follows. We start by setting $f_0(v) \defeq \tbl(v,0)$ for all $v \in V$. At the $n$-th iteration, we let $\Const_n \subseteq \Const$ be the set of all constraints violated by $f_n$. We then pick some $D$-independent\footnote{Recall that a set of vertices in a graph is \emphd{independent} if no two of its members are adjacent.} subset $\I_n \subseteq \Const_n$ 
    and, for each constraint $\const \in \I_n$ and every variable $v \in \dom(\const)$, update $f_n(v)$ to be the next value in the sequence $\tbl(v)$ (we use a function $\level_n \colon V \to \N$ to keep track of the current position in the sequence). This produces the labeling $f_{n+1}$. Taking a natural limit yields a (possibly partial) labeling $f$. The details are presented in Algorithm~\ref{alg:MT}. 

    {
    \floatname{algorithm}{Algorithm}
    \begin{algorithm}[t]\DontPrintSemicolon
        \caption{Moser--Tardos Algorithm ($\MTA$)}\label{alg:MT}
        \begin{flushleft}
            \textbf{Input:} A CSP $\Pi = (V, \Lambda, \Const, \dom, \Bad)$ and a table $\tbl \colon V \to \Lambda^\N$ \\ 
        \end{flushleft}

        \smallskip
        
             Initialize $\level_0(v) \defeq 0$ for all $v \in V$.\;

            \smallskip
            
            \For{$n =0$, $1$, $2$, \ldots}{

            \smallskip
            
                Define $f_{n}(v) \defeq \tbl(v, \level_{n}(v))$ for all $v \in V$.\;

                \smallskip

                 Let $\Const_n \defeq \set{\const \in \Const \,:\, \text{$f_n$ violates $\const$}}$ and pick a $D$-independent subset $\I_n \subseteq \Const_n$.\;

                \smallskip

                \For{$v \in V$}{
                
                \eIf{{\normalfont$v \in \dom(\const)$ for some $\const \in \I_n$}}{

                \smallskip

                 $\level_{n+1}(v) \defeq \level_n(v) + 1$\;}{

                 $\level_{n+1}(v) \defeq \level_n(v)$\;}
                }
            }

            \smallskip

             Let $\displaystyle \level(v) \defeq \lim_{n \to \infty} \level_n(v) \in \N \cup \set{\infty}$ and define $f(v) \defeq \tbl(v, \level(v))$ whenever $\level(v) < \infty$.

        \smallskip
    \end{algorithm}
    }

    Note that Algorithm~\ref{alg:MT} does not specify how the $D$-independent subsets $\I_n \subseteq \Const_n$ are to be chosen. Hence, the input of $\MTA$ does not determine the execution process of the algorithm uniquely. We call a sequence $\I = (\I_n)_{n \in \N}$ of $D$-independent subsets of $\Const$ an \emphd{$\MT$-sequence}, and say that an $\MT$-sequence is \emphd{consistent} with $(\Pi,\tbl)$ if it can be produced by $\MTA$ 
    with input $(\Pi, \tbl)$. 
    It is also possible to give a direct combinatorial description of $\MT$-sequences consistent with $(\Pi, \tbl)$. To this end, given an $\MT$-sequence $\I = (\I_n)_{n \in \N}$ and $\const \in \Const$, we write
    \[
        L_n(\I, \const) \,\defeq\, \set{i < n \,:\, \const \in \I_i}. 
    \]
    %

    \begin{Lemma}\label{lemma:tf}
        Let $\tbl$ be a table and let $\I$ be an $\MT$-sequence. For each $n \in \N$ and $v \in V$, let
        \begin{equation}\label{eq:tf}
            \level_n(\I, v) \,\defeq\, \sum_{\const \,\in\, \dom^{-1}(v)} |L_n(\I, \const)| 
            \qquad \text{and} \qquad f_n(\I, v) \,\defeq\, \tbl(v, \level_n(v)).
        \end{equation}
        Then $\I$ is consistent with $(\Pi,\tbl)$ if and only if for all $n \in \N$, $f_n(\I, \cdot)$ violates every constraint in $\I_n$.
    \end{Lemma}
    \begin{scproof}
        If we run $\MTA$ on input $(\Pi, \tbl)$ with $\I$ as the $\MT$-sequence, then the resulting functions $\level_n \colon V \to \N$ and $f_n \colon V \to \Lambda$ will be precisely given by \eqref{eq:tf}. Therefore, choosing $\I_n$ on step $n$ of the algorithm is allowed if and only if $\I_n \subseteq \Const_n$, i.e., if $f_n$ violates every constraint in $\I_n$. 
    \end{scproof}

    The \emphd{Maximal Moser--Tardos Algorithm}, or $\MMTA$ for short, is a variant of $\MTA$ with the additional requirement that each set $\I_n$ be an \emph{\ep{inclusion-}maximal} $D$-independent subset of $\Const_n$ \ep{this is \cite[Alg.~1.2]{MT}}. 
    We write $\MTA(\Pi, \tbl)$ (resp.~$\MMTA(\Pi, \tbl)$) for the set of all partial labelings $f$ that can be generated by $\MTA$ (resp.~$\MMTA$) on input $(\Pi, \tbl)$. 
    By definition, $\MMTA(\Pi, \tbl) \subseteq \MTA(\Pi, \tbl)$.

    \begin{Lemma}\label{lemma:limit}
        Suppose $\P[\Bad(\const)] < 1$ for all $\const \in \Const$ and let $\tbl$ be a table. 
        If a labeling $f \in \MMTA(\Pi, \tbl)$ is defined on all of $V$, then $f$ is a solution to $\Pi$.
    \end{Lemma}
    \begin{scproof}
        Note that if the domain of a constraint $\const$ is empty, then there are two options for $\Bad(\const)$, namely $\0$ and $\set{\0}$, and we have $\P[\0] = 0$ and $\P[\set{\0}] = 1$. Since $\P[\Bad(\const)] < 1$ by assumption, we conclude that $\Bad(\const) = \0$. In other words, such a constraint $\const$ is satisfied by every labeling.
        
        Now suppose $f \in \MMTA(\Pi, \tbl)$ is defined everywhere and violates some constraint $\const \in \Const$. It follows that $\dom(\const) \neq \0$. Let $\level_n$ and $\level$ be as in Algorithm~\ref{alg:MT}, and let $n \in \N$ be such that $\level(v) = \level_n(v)$ for all $v \in \dom(\mathfrak{c})$, which exists because $\level$ is finite everywhere. Then $f$ agrees with $f_n$ on $\dom(\const)$, so $f_n$ also violates $\const$. This implies that $\const \in \Const_n$, and, by the maximality of $\I_n$ and since $\dom(\const) \neq \0$, there must be some $\const' \in \I_n$ with $\const' \in N[\const]$. 
        But then for any $v \in \dom(\const) \cap \dom(\const')$, we have $\level_{n+1}(v) = \level_n(v) + 1$, contradicting the choice of $n$.
    \end{scproof}


    \subsection{Good and locally good tables}\label{subsec:good_definitions}

    It follows from Lemma~\ref{lemma:limit} that to solve $\Pi$, it is enough to find a table $\tbl$ such that some labeling $f \in \MMTA(\Pi, \tbl)$ is defined everywhere. We shall seek a table with the following stronger property:

    \begin{defn}[Good tables]
        A table $\tbl$ is \emphd{good} if every $f \in \MTA(\Pi, \tbl)$ is defined on all of $V$.
    \end{defn}
    

    The following is the main result of  \cite{MT}:

    \begin{theo}[{Moser--Tardos \cite[Thm.~1.2]{MT}}]\label{theo:MT}
        Let $p$ and $d$ be as in Theorem~\ref{theo:LLL} \ep{i.e., the LLL}. 
        Sample a random table $\tbl \colon V \to \Lambda^\N$ from the product space $(\Lambda^{V \times \N}, \P^{V \times \N})$. Then for each $v \in V$, \[\P[\text{every $f \in \MTA(\Pi, \tbl)$ is defined on $v$}] \,=\, 1.\]
    \end{theo}

    It follows that if $V$ is finite \ep{or countable}, then a random table $\tbl \colon V \to \Lambda^\N$ is almost surely good; in particular, a good table exists. As explained in Remark~\ref{remk:LLL}\ref{item:compactness}, Theorem~\ref{theo:LLL} for infinite $V$ then follows via a compactness argument. Moreover, although the proof of Theorem~\ref{theo:MT} given in \cite{MT} assumes $V$ is finite, it works just as well in the infinite case; see \cite[\S3]{bernshteyn2019measurable} for a detailed presentation of the argument for arbitrary $V$. (This was also observed in \cite{RSh,Kun}.)

    Unfortunately, even though Theorem~\ref{theo:MT} holds for infinite $V$, it cannot be applied directly to find a Borel solution to $\Pi$, for two reasons. First, if $V$ is uncountable, we can no longer conclude that a random table is good with probability $1$ (or even with positive probability). Second---and more importantly---a random table would typically not be a Borel function. As mentioned in the introduction, Cs\'oka \emph{et al.}~\cite{CGMPT} handled this issue by using the same random values for multiple variables. Namely, they partitioned $V$ into finitely many ``sparse'' Borel subsets as $V = V_1 \cup \ldots \cup V_k$ and generated a table $\tbl \colon V \to \Lambda^\N$ by independently  picking $k$ random sequences $\tbl_1$, \ldots, $\tbl_k \in \Lambda^\N$ and setting $\tbl(v) \defeq \tbl_i$ for all $v \in V_i$. They then were able to establish a suitable analog of Theorem~\ref{theo:MT} for this construction under the assumptions of Theorem~\ref{theo:CGMPT_slow}.
    
    Our approach instead is to treat the problem of finding a good table $\tbl \colon V \to \Lambda^\N$ as \emph{a CSP in its own right}. The idea is to ``localize'' the notion of goodness by ``zooming in'' on $r$-balls around individual constraints in the graph $D$. To this end, given $\const \in \Const$ and $r \in \N$, we define the \emphd{$(\const, r)$-local CSP} $\Pi_{\const, r} \defeq (V, \Lambda, \Const, \dom, \Bad_{\const, r})$ by letting, for all $\mathfrak{a} \in \Const$,
    \[
        \Bad_{\const, r}(\mathfrak{a}) \,\defeq\, \begin{cases}
            \Bad(\mathfrak{a}) &\text{if } \mathfrak{a} \in B_D(\const, r),\\
            \Lambda^{\dom(\mathfrak{a})} &\text{otherwise}.
        \end{cases}
    \]
    That is, from the point of view of $\Pi_{\const, r}$, the constraints at distance more than $r$ from $\const$ in $D$ are always violated. In particular, to check whether an $\MT$-sequence is consistent with $(\Pi_{\const, r}, \tbl)$, we only need to consider the constraints in $B_D(\const, r)$:

    \begin{Lemma}\label{lemma:tf_local}
        Let $\tbl$ be a table and let $\I$ be an $\MT$-sequence. For each $n \in \N$ and $v \in V$, define $\level_n(\I, v)$ and $f_n(\I,v)$ as in \eqref{eq:tf}. Then, for any $\const \in \Const$ and $r \in \N$,  $\I$ is consistent with $(\Pi_{\const, r},\tbl)$ if and only if for all $n \in \N$, $f_n(\I, \cdot)$ violates every constraint in $\I_n \cap B_D(\const, r)$.
    \end{Lemma}
    \begin{scproof}
        Follows immediately by applying Lemma~\ref{lemma:tf} with $\Pi_{\const,r}$ in place of $\Pi$.
    \end{scproof}

    It is clear that, unless $B_D(\const, r) = \Const$, the CSP $\Pi_{\const, r}$ has no solutions, because the constraints in $\Const \setminus B_D(\const, r)$ are, by definition, always violated. In particular, the Moser--Tardos Algorithm must fail to generate a solution to $\Pi_{\const,r}$. However, we may ask whether the constraints in $\Const \setminus B_D(\const, r)$ are in some sense the {``main reason''} for this failure. Namely, if $\I = (\I_n)_{n \in \N}$ is an $\MT$-sequence consistent with $(\Pi_{\const, r}, \tbl)$, must it be that ``many'' constraints in the sets $\I_n$ come from outside $B_D(\const, r)$?
    
    In order to quantify the word ``many'' in the preceding sentence, we restrict our attention to \emphd{finite} $\MT$-sequences, i.e., $\MT$-sequences $\I = (\I_n)_{n \in \N}$ such that $\sum_{n \in \N} |\I_n| < \infty$. In other words, $\I$ is finite if each $\I_n$ is a finite set and $\I_n = \0$ for all large enough $n \in \N$. 
    Now, if $\I = (\I_n)_{n \in \N}$ is a finite $\MT$-sequence consistent with $(\Pi_{\const, r}, \tbl)$, we can consider the sets
    \[
        \I_n \cap B_D(\const, r) \qquad  \text{and} \qquad \I_n \setminus B_D(\const,r)
    \]
    for each $n \in \N$, and we want the former ones to not be much larger than the latter. 
    %
    Formally, we say a finite $\MT$-sequence $\I = (\I_n)_{n \in \N}$ is \emphd{$(\const, r, N, \epsilon)$-F\o{}lner} for $\const \in \Const$, $r$, $N \in \N$, and $\epsilon \in (0,1)$ if 
    \[
        \sum_{n \in \N} |\I_n \cap B_D(\const, r)| \,\geq\, N \qquad \text{and} \qquad  \sum_{n \in \N} |\I_n \setminus B_D(\const, r)| \,<\, \epsilon \sum_{n \in \N}  |\I_n|. 
    \]

    \begin{defn}[Locally good tables]
        Fix $\const \in \Const$, $R$, $N \in \N$, and $\epsilon \in (0,1)$. A table $\tbl \colon V \to \Lambda^\N$ is \emphd{$(\const, R, N, \epsilon)$-locally good} if for all $0 \leq r < R$, there is no 
        $(\const, r, N, \epsilon)$-F\o{}lner $\MT$-sequence 
        consistent with $(\Pi_{\const, r}, \tbl)$. 
        A table $\tbl$ is \emphd{$(R,N,\epsilon)$-locally good} if it is $(\const, R, N, \epsilon)$-locally good for all $\const \in \Const$.
    \end{defn}

    The following lemma, proved in \S\ref{subsec:good}, is the only place where we use a bound on $\egr(D)$:

    \begin{Lemma}\label{lemma:good}
        Let $R \in \N$ and $\epsilon \in (0,1)$ be such that $\gamma_D(R) < (1-\epsilon)^{-R}$, where $\gamma_D$ is the growth function of $D$. If a table $\tbl \colon V \to \Lambda^\N$ is $(R,N,\epsilon)$-locally good for some $N \in \N$, then $\tbl$ is good.
    \end{Lemma}

    Next we construct a CSP $\mathsf{LG}(R,N,\epsilon)$ with label set $\Lambda^\N$ such that a table $\tbl \colon V \to \Lambda^\N$ is a solution to $\mathsf{LG}(R,N,\epsilon)$ if and only if $\tbl$ is $(R,N,\epsilon)$-locally good. To this end, we define 
    \[
        \mathsf{LG}(R,N,\epsilon) \,\defeq\, \big(V, \, \Lambda^\N, \, \Const, \, \dom_R, \, \LBad_{R,N,\epsilon}\big),
    \]
    where for each $\const \in \Const$, $\dom_R(\const) \defeq \bigcup_{\mathfrak{a} \in B_D(\const, R)} \dom(\mathfrak{a})$ and $\LBad_{R,N,\epsilon}$ is given by the following:

    \begin{Lemma}\label{lemma:local_CSP}
        Fix $\const \in \Const$, $R$, $N \in \N$, and $\epsilon \in (0,1)$. If tables $\tbl_0$, $\tbl_1 \colon V \to \Lambda^\N$ agree on $\dom_R(\const)$, then
        \[
            \text{$\tbl_0$ is $(\const, R,N,\epsilon)$-locally good} \quad \Longleftrightarrow \quad \text{$\tbl_1$ is $(\const, R,N,\epsilon)$-locally good}.
        \]
        Therefore, there is a set $\LBad_{R, N, \epsilon}(\const) \subseteq (\Lambda^\N)^{\dom_R(\const)}$ such that for any $\phi \colon \dom_R(\const) \to \Lambda^\N$, the following statements are equivalent:
        \begin{enumerate}[label=\ep{\normalfont\roman*}]
            \item some table extending $\phi$ is not $(\const, R, N, \epsilon)$-locally good,
            \item every table extending $\phi$ is not $(\const, R, N, \epsilon)$-locally good,
            \item $\phi \in \LBad_{R, N, \epsilon}(\const)$.
        \end{enumerate}
    \end{Lemma}
    \begin{scproof}
        If $\tbl_0$ is not $(\const, R, N, \epsilon)$-locally good, then for some $r < R$, there is a $(\const, r, N, \epsilon)$-F\o{}lner $\MT$-sequence $\I$ consistent with $(\Pi_{\const, r}, \tbl_0)$. Since $\tbl_0$ and $\tbl_1$ agree on $\dom_r(\const)$, Lemma~\ref{lemma:tf_local} shows that $\I$ is consistent with $(\Pi_{\const, r}, \tbl_1)$ as well. Hence, $\tbl_1$ is also not $(\const, R, N, \epsilon)$-locally good, as claimed.
    \end{scproof}

    It follows from Lemma~\ref{lemma:local_CSP} that a table $\tbl \colon V \to \Lambda^\N$ is a solution to $\mathsf{LG}(R, N, \epsilon)$ if and only if $\tbl$ is $(R,N,\epsilon)$-locally good, as intended.
    
    Next we proceed to bound the maximum degree of the dependency graph of $\mathsf{LG}(R, N, \epsilon)$: 

    \begin{Lemma}\label{lemma:LG_degree}
        For any $R$, $N \in \N$ and $\epsilon > 0$, the maximum degree of the dependency graph of $\mathsf{LG}(R, N, \epsilon)$ is at most
        $\gamma_D(2R+1) - 1$. 
    \end{Lemma}
    \begin{scproof}
        Distinct constraints $\const$, $\const' \in \Const$ are adjacent in the dependency graph of $\mathsf{LG}(R, N, \epsilon)$ if and only if $\dom_R(\const) \cap \dom_R(\const') \neq \0$, which means that for some $\mathfrak{a} \in B_D(\const,R)$ and $\mathfrak{a}' \in B_D(\const',R)$, we have $\dom(\mathfrak{a}) \cap \dom(\mathfrak{a}') \neq \0$, i.e., $\mathfrak{a}' \in B_D(\mathfrak{a},1)$. The latter condition implies that $\const' \in B_D(\const, 2R+1) \setminus \set{\const}$. Since $|B_D(\const, 2R+1)| \leq \gamma_D(2R+1)$, the result follows.
    \end{scproof}

    The only consequence of Lemma~\ref{lemma:LG_degree} we need is that the maximum degree of the dependency graph of $\mathsf{LG}(R, N, \epsilon)$ can be bounded by a function of $R$ \emph{independent of $N$}. By contrast, we will show that $\P[\LBad_{R, N, \epsilon}(\const)]$ is bounded above by a quantity that \emph{goes to $0$ as $N \to \infty$} (in fact, it goes to $0$ exponentially quickly, but that is not important for our purposes):

    \begin{Lemma}\label{lemma:LG_prob}
        Suppose there exist $p \in [0,1)$, $d \in \N$, $s > 1$, and $\epsilon$, $\eta \in (0,1)$ such that:
        \begin{enumerate}[label=\ep{\normalfont\roman*}]
            \item $\P[\Bad(\const)] \leq p$ for all $\const \in \Const$, the maximum degree of $D$ is at most $d$, and  $p\, (\e \,(d+1))^{s} < 1$, 
            \item $\displaystyle \epsilon + \frac{1}{s} < 1$ and $\displaystyle p^{1 - \epsilon - \frac{1}{s}} \leq \frac{1-\eta}{1+\eta}$.
        \end{enumerate}
        Then for all $\const \in \Const$ and $R$, $N \in \N$, we have
        $\P[\LBad_{R,N,\epsilon}(\const)] \leq (1+\eta)^{-N} \,F(d, \eta, R)$, where the quantity $F(d, \eta, R)$ is independent of $N$.
    \end{Lemma}

    We prove Lemma~\ref{lemma:LG_prob} in \S\ref{subsec:witness}.

     Using Lemmas~\ref{lemma:LG_degree} and \ref{lemma:LG_prob}, for any fixed $R$, we are able to pick $N \in \N$ so large that $\mathsf{LG}(R,N,\epsilon)$ satisfies the assumptions of Lemma~\ref{lemma:double_exp}. As a result, we can find a Borel $(R,N,\epsilon)$-locally good table $\tbl$---which must be good by Lemma~\ref{lemma:good}---and use it to finish the proof of Theorem~\ref{theo:main}. The details of this part of the argument, as well as the proof of Lemma~\ref{lemma:double_exp}, are presented in \S\ref{sec:Borel}.


    \subsection{Preliminary facts about tables and $\MT$-sequences}\label{subsec:facts}

    The following observation will be useful:

    \begin{Lemma}\label{lemma:restrict}
        Fix $\const \in \Const$ and $r \in \N$. Let $\I = (\I_n)_{n \in \N}$ be an $\MT$-sequence consistent with $(\Pi_{\const, r}, \tbl)$ for some table $\tbl$. Then, for any $R > r$, $\I' \defeq (\I_n \cap B_D(\const, R))_{n \in \N}$ is also consistent with $(\Pi_{\const, r}, \tbl)$.
    \end{Lemma}
    \begin{scproof}
        Recall for all $\mathfrak{a} \in \Const$ and $n \in \N$ the notation
        \[
            L_n(\I, \mathfrak{a}) \,\defeq\, \set{i < n \,:\, \mathfrak{a} \in \I_i} \qquad \text{and} \qquad L_n(\I', \mathfrak{a}) \,\defeq\, \set{i < n \,:\, \mathfrak{a} \in \I_i'}.
        \]
        For each $n \in \N$ and $v \in V$, let $f_n(\I, v)$ and $f_n(\I', v)$ be defined according to \eqref{eq:tf}. By Lemma~\ref{lemma:tf_local}, we need to show that for all $\mathfrak{a} \in \I'_n \cap B_D(\const, r) = \I_n \cap B_D(\const, r)$, the function $f_n(\I', \cdot)$ violates $\mathfrak{a}$. To this end, note that if $v \in \dom(\mathfrak{a})$, then every constraint $\mathfrak{b} \in \dom^{-1}(v)$ belongs to $B_D(\const, r+1) \subseteq B_D(\const, R)$, and hence $L_n(\I', \mathfrak{b}) = L_n(\I, \mathfrak{b})$. It follows that  $f_n(\I', v) = f_n(\I, v)$. Since $\I$ is consistent with $(\Pi_{\const, r}, \tbl)$, $f_n(\I, \cdot)$ violates $\mathfrak{a}$ by Lemma~\ref{lemma:tf_local}, which implies that $f_n(\I', \cdot)$ violates $\mathfrak{a}$ as well, as desired.
    \end{scproof}

    A consequence of Lemma~\ref{lemma:restrict} is that the property of being $(\const, R, N, \epsilon)$-locally good for a table $\tbl$ is ``localized'' in the graph $D$, in the sense that it only involves the constraints in $B_D(\const, R)$. To make this statement precise, we say that an $\MT$-sequence $\I = (\I_n)_{n \in \N}$ is \emphd{$(\const, R)$-bounded} for some $\const \in \Const$ and $R \in \N$ if $\I_n \subseteq B_D(\const, R)$ for all $n \in \N$.

    \begin{Lemma}\label{lemma:restricted_Folner}
        Fix $\const \in \Const$, $R$, $N \in \N$, and $\epsilon \in (0,1)$. A table $\tbl$ is $(\const, R, N, \epsilon)$-locally good if and only if for all $r < R$, no $(\const, R)$-bounded $(\const, r, N, \epsilon)$-F\o{}lner $\MT$-sequence $\I$ is consistent with $(\Pi_{\const, r}, \tbl)$.
    \end{Lemma}
    \begin{scproof}
        If $\tbl$ is $(\const, R, N, \epsilon)$-locally good, then such an $\MT$-sequence $\I$ does not exist by definition. Now suppose $\tbl$ is not $(\const, R,N,\epsilon)$-locally good and let $\I = (\I_n)_{n \in \N}$ be a $(\const, r, N, \epsilon)$-F\o{}lner $\MT$-sequence consistent with $(\Pi_{\const, r}, \tbl)$. Then the $\MT$-sequence $\I' \defeq (\I_n \cap B_D(\const, R))_{n \in \N}$ is also consistent with $(\Pi_{\const, r}, \tbl)$ by Lemma~\ref{lemma:restrict}, and it is $(\const, r, N, \epsilon)$-F\o{}lner and $(\const, R)$-bounded, as desired.
    \end{scproof}

    \subsection{Locally good implies good: Proof of Lemma~\ref{lemma:good}}\label{subsec:good}

    Let $R \in \N$ and $\epsilon \in (0,1)$ be such that $\gamma_D(R) < (1-\epsilon)^{-R}$ and let $\tbl \colon V \to \Lambda^\N$ be an $(R,N,\epsilon)$-locally good table for some $N \in \N$. Our goal is to show that $\tbl$ is good. 
    
    Assume, toward a contradiction, that there is a labeling $f \in \MTA(\Pi, \tbl)$ whose domain is not all of $V$ and let $\I = (\I_n)_{n \in \N}$ be the $\MT$-sequence consistent with $(\Pi, \tbl)$ used to generate $f$. Recall for all $\const \in \Const$ and $n \in \N$ the notation $L_n(\I, \const) \,\defeq\, \set{i < n \,:\, \const \in \I_i}$.  If $\level_n$ and $\level$ are as in Algorithm~\ref{alg:MT}, then 
    \[
        \level_n(v) \,=\, \sum_{\const \,\in\, \dom^{-1}(v)} |L_n(\I, \const)|.
    \]
    Since $\level(v) = \infty$ for some $v \in V$ and $|\dom^{-1}(v)| \leq d+1$ for all $v \in V$, there must be $n \in \N$ and $\const \in \Const$ such that
    $
        |L_n(\I, \const)| > N
    $.
    We fix any such $n \in \N$ and choose $\const \in \Const$ to maximize $|L_n(\I, \const)|$.

    For each integer $-1 \leq r < R$, we define a finite $\MT$-sequence $\I^r = (\I^r_i)_{i \in \N}$ via
    \[
        \I^r_i \,\defeq\, \begin{cases}
            \I_i \cap B_D(\const, r+1) &\text{if } i < n,\\
            \0 &\text{if } i \geq n.
        \end{cases}
    \]
    Consider any $0 \leq r < R$. Since $\I$ is consistent with $(\Pi, \tbl)$ and hence also with $(\Pi_{\const, r}, \tbl)$, Lemma~\ref{lemma:restrict} shows that $\I^r$ is consistent with $(\Pi_{\const,r}, \tbl)$. Furthermore, note that
    \[
        \sum_{i \in \N} |\I_i^r \cap B_D(\const, r)| \,=\, \sum_{i \in \N}  |\I_i^{r-1}| \,\geq\, \sum_{i \in \N} |\I^{-1}_i| \,=\, |L_n(\I, \const)| \,>\, N.
    \]
    Since $\tbl$ is $(R,N,\epsilon)$-locally good, $\I^r$ cannot be $(\const, r, N, \epsilon)$-F\o{}lner, and hence 
    \[
        \sum_{i \in \N} |\I_i^r \setminus B_D(\const, r)| \,\geq\, \epsilon \sum_{i \in \N}  |\I_i^r|.\]
    This implies that, for $0 \leq r < R$,
    \[
        \sum_{i \in \N} |\I_i^r| \,\geq\, \frac{1}{1-\epsilon} \, \sum_{i \in \N}  |\I_i^r \cap B_D(\const, r)| \,=\, \frac{1}{1-\epsilon} \, \sum_{i \in \N}  |\I_i^{r-1}|.
    \]
    Applying this inequality iteratively, we obtain 
    \begin{equation}\label{eq:IR_lower}
        \sum_{i \in \N} |\I_i^{R-1}| \,\geq\, \left(\frac{1}{1 - \epsilon}\right)^{R} \sum_{i \in \N} |\I^{-1}_i| \,=\, \left(\frac{1}{1 - \epsilon}\right)^{R} |L_n(\I, \const)|.
    \end{equation}
    On the other hand,
    \begin{equation}\label{eq:IR_upper}
        \sum_{i \in \N} |\I_i^{R-1}| \,=\, \sum_{\mathfrak{a} \in B_D(\const, R)} |L_n(\I, \mathfrak{a})| \,\leq\, |B_D(\const, R)| \, |L_n(\I, \const)| \,\leq\, \gamma_D(R) \, |L_n(\I, \const)|
    \end{equation}
    by the choice of $\const$. Combining \eqref{eq:IR_lower} and \eqref{eq:IR_upper}, we get $\gamma_D(R) \geq (1-\epsilon)^{-R}$, which is a contradiction. The proof of Lemma~\ref{lemma:good} is complete.

    \subsection{Witness digraphs and the proof of Lemma~\ref{lemma:LG_prob}}\label{subsec:witness}
    
    Our proof of Lemma~\ref{lemma:LG_prob} relies, in a black box manner, on certain basic results originating in the work of Moser and Tardos \cite{MT}. The central idea of Moser--Tardos theory is to analyze $\MTA$ by means of certain \emph{witness structures} that keep track of $\MTA$'s execution process. One can find several 
    types of such structures in the literature, e.g., the \emph{witness trees} from the original paper \cite{MT}, \emph{stable set sequences} \cite{KolipakaSzegedy}, \emph{neat piles} \cite{bernshteyn2019measurable}, and \emph{witness digraphs} \cite{haeupler2017parallel}. 
    The latter will be especially convenient for our purposes, although it must be stressed that all these different structures are essentially equivalent and can be straightforwardly converted into each other.
    
    The definition of a witness digraph appeared implicitly in the paper \cite{KolipakaSzegedy} by Kolipaka and Szegedy and was formalized in \cite{haeupler2017parallel} by Haeupler and Harris. 
    See \cite{harris2020new, harris2021new, harris2023deterministic, he2023moser}
    for a selection of other works that rely on this notion. 

    In the sequel, a \emphd{digraph} means a simple directed graph. A \emphd{decorated digraph} is a digraph $\Gamma$ equipped with a mapping $\updelta \colon V(\Gamma) \to \Const$, called a \emphd{decoration}. We shall always use the symbol $\updelta$ to denote the decoration on a decorated digraph. 
    Isomorphisms between decorated digraphs are required to preserve the decoration. 

    To each finite $\MT$-sequence we associate a decorated digraph called its \emph{full witness digraph}:

    \begin{defn}[Full witness digraph]
        Let $\I = (\I_n)_{n \in \N}$ be a finite $\MT$-sequence. The \emphd{full witness digraph} of $\I$ is the finite decorated digraph $\Gamma(\I)$ defined as follows. The vertex set of $\Gamma(\I)$ is 
        \[
            V(\Gamma({\I})) \,\defeq\, \big\{(n, \const) \,:\, n \in \N, \, \const \in \I_n\big\}.
        \]
        There is a directed edge from $(n, \const)$ to $(n', \const')$ in $\Gamma({\I})$ if and only if $n < n'$ and $\const \in N[\const']$. The decoration $\updelta$ of $\Gamma({\I})$ is given by $\updelta(n, \const) \defeq \const$ for all vertices $(n, \const)$.
    \end{defn}


    \begin{Lemma}\label{lemma:witness}
        Let $\I = (\I_n)_{n \in \N}$ be a finite $\MT$-sequence and let $\Gamma \defeq \Gamma(\I)$. Then:
        \begin{enumerate}[label=\ep{\normalfont\roman*}]
            \item\label{item:acyclic} $\Gamma$ has no directed cycles,
            \item\label{item:cliques} for all $x$, $y \in V(\Gamma)$, $(x,y) \in E(\Gamma)$ or $(y,x) \in E(\Gamma)$ if and only if $x \neq y$ and $\updelta(x) \in N[\updelta(y)]$. 
        \end{enumerate}
    \end{Lemma}
    \begin{scproof}
        If $(x,y) \in E(\Gamma)$, then the first coordinate of $x$ is strictly less than the first coordinate of $y$, so no directed cycle is possible. Item \ref{item:cliques} is immediate from the way the edge set of $\Gamma$ is defined.
    \end{scproof}

    We call a finite decorated digraph $\Gamma$ satisfying conditions \ref{item:acyclic} and \ref{item:cliques} of Lemma~\ref{lemma:witness} a \emphd{witness digraph}. Thus, Lemma~\ref{lemma:witness} says that the full witness digraph of a finite $\MT$-sequence is a witness digraph. 
    Now we relate witness digraphs and tables:

    \begin{defn}[Compatible witness digraphs]
        A witness digraph $\Gamma$ is \emphd{compatible} with $(\Pi, \tbl)$, where $\tbl$ is a table, if there exists an $\MT$-sequence $\I$ consistent with $(\Pi,\tbl)$ such that $\Gamma \cong \Gamma(\I)$. 
    \end{defn}



    
    The next result, which we call the \emph{First Fundamental Theorem} of Moser--Tardos theory, computes the probability that a given witness digraph is compatible with $(\Pi, \tbl)$ for random $\tbl$:

    \begin{theo}[{First Fundamental Theorem \cite[Prop.~2.5]{haeupler2017parallel}}]\label{theo:MT1}
        Sample a random table $\tbl \colon V \to \Lambda^\N$ from the space $(\Lambda^{V \times \N}, \P^{V \times \N})$. Then, for every witness digraph $\Gamma$, we have
        \[
            \P[\Gamma \text{ is compatible with } (\Pi,\tbl)] \,=\, \prod_{x \in V(\Gamma)} \P[\Bad(\updelta(x))].
        \]
    \end{theo}

%
        Given a constraint $\const \in \Const$, we let $\Sink(\const)$ be the class of all witness digraphs $\Gamma$ such that $\Gamma$ has a single sink $\sigma \in V(\Gamma)$ and $\updelta(\sigma) = \const$. 
        We also let $\Sink^\ast(\const)$ be a set 
        containing exactly one representative from each isomorphism class of witness digraphs in $\Sink(\const)$. 
    The \emph{Second Fundamental Theorem} of Moser--Tardos theory is a purely algebraic statement that allows one to bound the quantities appearing in Theorem~\ref{theo:MT1} for $\Gamma \in \Sink^*(\const)$:

    \begin{theo}[{Second Fundamental Theorem \cite[Prop.~2.10]{haeupler2017parallel}}]\label{theo:MT2}
        Suppose $\alpha$, $\beta \colon \Const \to [0,1)$ are functions such that for all $\const \in \Const$,
        \[
            \alpha(\const) \,\leq\, \beta(\const) \prod_{\const' \in N(\const)} (1 - \beta(\const')).
        \]
        Then for each $\const \in \Const$,
        \[
            \sum_{\Gamma \in \Sink^\ast(\const)}  \, \prod_{x \in V(\Gamma)} \alpha(\updelta(x)) \,\leq\, \frac{\beta(\const)}{1 - \beta(\const)}.
        \]
    \end{theo}

    Theorems~\ref{theo:MT1} and \ref{theo:MT2} are essentially present in the seminal paper \cite{MT}, cf.~\cite[Lem.~2.1 and 3.1]{MT}; the formulations using witness digraphs were given in \cite{haeupler2017parallel}.

    We now have all the tools needed to prove Lemma~\ref{lemma:LG_prob}.

    \begin{scproof}[ of Lemma~\ref{lemma:LG_prob}]
        Fix $p \in [0,1)$, $d \in \N$, $s > 1$, and $\epsilon$, $\eta \in (0,1)$ as in the statement of Lemma~\ref{lemma:LG_prob}.  
        Let $\const \in \Const$ and $R$, $N \in \N$. Our goal is to bound 
        $\P[\LBad_{R,N,\epsilon}(\const)]$, i.e., the probability that a random table $\tbl \colon V \to \Lambda^\N$ is not $(\const, R, N, \epsilon)$-locally good. We shall achieve this by applying Theorem~\ref{theo:MT2} to certain auxiliary CSPs $\Sigma_{r}$ in place of $\Pi$. Specifically, for each $r < R$, we define $\Sigma_{r}$ by adding to $\Pi_{\const, r}$ a new constraint $\mathfrak{u}$ with $\dom(\mathfrak{u}) \defeq \dom_{R-1}(\const)$ that is always violated.

        \begin{claim}\label{claim:found_Gamma}
            If a table $\tbl$ is not $(\const, R, N, \epsilon)$-locally good, then for some $r < R$, there is a witness digraph $\Gamma \in \Sink^*(\mathfrak{u})$ compatible with $(\Sigma_r, \tbl)$ such that 
            \begin{itemize}
                \item $|\set{x \in V(\Gamma) \,:\, \updelta(x) = \mathfrak{u}}| = 1$,
                \item $|\set{x \in V(\Gamma) \,:\, \updelta(x) \in B_D(\const, r)}| \geq N$, and
                \item $|\set{x \in V(\Gamma) \,:\, \updelta(x) \in \Const \setminus B_D(\const, r)}| < \epsilon (|V(\Gamma)| - 1)$.
            \end{itemize}
        \end{claim}
        \begin{claimproof}
            By Lemma~\ref{lemma:restricted_Folner}, there exist $r < R$ and a $(\const,R)$-bounded $(\const, r, N, \epsilon)$-F\o{}lner $\MT$-sequence $\I = (\I_n)_{n \in \N}$ consistent with $(\Pi_{\const, r}, \tau)$. Fix $k \in \N$ such that $\I_n = \0$ for all $n \geq k$ and define
        \[
            \I^*_n \,\defeq \, \begin{cases}
                \I_n &\text{if } n < k,\\
                \set{\mathfrak{u}} &\text{if } n = k, \\
                \0 &\text{if } n > k.
            \end{cases}
        \]
        This gives a finite $\MT$-sequence $\I^* \defeq (\I^*_n)_{n \in \N}$ for $\Sigma_r$. Observe that $\I^*$ is consistent with $(\Sigma_r, \tbl)$, because $\I$ is consistent with $(\Pi_{\const, r}, \tbl)$ and $\mathfrak{u}$ is always violated. Now we take $\Gamma \cong \Gamma(\I^*)$, noting that the vertex $(k, \mathfrak{u})$ is the unique sink in $\Gamma(\I^*)$. 
        \end{claimproof}

        Let $\mathbb{G}_r \subseteq \Sink^*(\mathfrak{u})$ be the set of all witness digraphs $\Gamma$ satisfying the properties in Claim~\ref{claim:found_Gamma}. It follows that, for a random table $\tbl$ sampled from $(\Lambda^{V \times \N}, \P^{V \times \N})$, 
        \[
            \P[\text{$\tbl$ is not $(\const, R, N, \epsilon)$-locally good}] \,\leq\, \sum_{r < R} \,\sum_{\Gamma \in \mathbb{G}_r} \P[\text{$\Gamma$ is compatible with $(\Sigma_r, \tbl)$}].
        \]
        Recall that, by the definition of $\Pi_{\const, r}$, the constraints in $\Const \setminus B_D(\const, r)$ are always violated, and the same goes for $\mathfrak{u}$. Since $\P[\B(\mathfrak{a})] \leq p$ for all $\mathfrak{a} \in B_D(\const, r)$, Theorem~\ref{theo:MT1} shows that for each $\Gamma \in \mathbb{G}_r$,
        \begin{align*}
            \P[\text{$\Gamma$ is compatible with $(\Sigma_r, \tbl)$}] \,&=\, \prod_{x \in V(\Gamma) \,:\, \updelta(x) \in B_D(\const, r)} \P[\Bad(\updelta(x))]\,\leq\, p^{|\set{x \in V(\Gamma) \,:\, \updelta(x) \in B_D(\const, r)}|}.
        \end{align*}
        By the definition of $\mathbb{G}_r$, the last quantity is less than $p^{(1-\epsilon)(|V(\Gamma)| - 1)}$. It remains to bound
        \[
            \sum_{r < R} \, \sum_{\Gamma \in \mathbb{G}_r} p^{(1-\epsilon)(|V(\Gamma)| - 1)}.
        \]
        To this end, we shall invoke Theorem~\ref{theo:MT2}.

        Define functions $\alpha_r$, $\beta_r \colon \Const \cup \set{\mathfrak{u}} \to [0,1)$ as follows. Note that there is $\zeta = \zeta(d) \in (0,1)$ such that \[\zeta \, (1-\zeta)^{d} \,\geq\, (\e\,(d+1))^{-1},\] namely $\zeta = 1/\e$ for $d= 0$ and $\zeta = 1/(d+1)$ for $d > 0$. We set 
        \[
            \xi \,\defeq\, \eta\, \left(1 - \zeta\right)^{\gamma_D(R)},
        \]
        and define
        \[
            \alpha_r(\mathfrak{a}) \,\defeq\, \begin{cases}
                (1+\eta)\,p^{1-\epsilon} &\text{if } \mathfrak{a} \in \Const,\\
                \xi &\text{if } \mathfrak{a} = \mathfrak{u}
            \end{cases} \qquad \text{and} \qquad \beta_r(\mathfrak{a}) \,\defeq\, \begin{cases}
                \zeta &\text{if } \mathfrak{a} \in \Const,\\
                \eta &\text{if } \mathfrak{a} = \mathfrak{u}.
            \end{cases}
        \]
        Using that $|V(\Gamma)| \geq N + 1$ for all $\Gamma \in \mathbb{G}_r$,  we can write
        \begin{align*}
            \sum_{\Gamma \in \mathbb{G}_r} p^{(1-\epsilon)(|V(\Gamma)| - 1)} \,&=\, \sum_{\Gamma \in \mathbb{G}_r}\frac{1}{\xi \, (1+\eta)^{|V(\Gamma)| - 1}} \, \prod_{x \in V(\Gamma)} \alpha_r(\updelta(x)) \\
            &\leq\, \frac{1}{\xi \, (1+\eta)^{N}} \, \sum_{\Gamma \in \Sink^*(\mathfrak{u})} \, \prod_{x \in V(\Gamma)} \alpha_r(\updelta(x)).
        \end{align*}
        Now we verify the assumptions of Theorem~\ref{theo:MT2}. In the dependency graph of $\Sigma_r$, the neighborhood of each $\mathfrak{a} \in \Const$ is either $N(\mathfrak{a})$ or $N(\mathfrak{a}) \cup \set{\mathfrak{u}}$, and
        \begin{align*}
            &\beta_r(\mathfrak{a}) \left(\prod_{\mathfrak{b} \in N(\mathfrak{a})} (1 - \beta_r(\mathfrak{b}))\right) (1 - \beta_r(\mathfrak{u})) \,\geq\, \zeta \left(1 - \zeta\right)^{d} (1 - \eta) \\
            &\hspace{3cm}\geq\, \frac{1 - \eta}{ \e \, (d+1)} \,>\, (1-\eta) p^{1/s} \,\geq\, (1+\eta)p^{1-\epsilon} \,=\, \alpha_r(\mathfrak{a}),
        \end{align*}
        where we use the bounds $p\,(\e\,(d+1))^s < 1$ and $p^{1 - \epsilon - \frac{1}{s}} \leq (1-\eta)/(1+\eta)$. 
        On the other hand, the neighborhood of $\mathfrak{u}$ is $B_D(\const, R)$, and
        \[
            \beta_r(\mathfrak{u})\prod_{\mathfrak{a} \in B_D(\const, R)}(1 - \beta_r(\mathfrak{a})) \,=\, \eta \, \left(1 - \zeta\right)^{|B_D(\const, R)|} \,\geq\, \xi \,=\, \alpha_r(\mathfrak{u}).
        \]
        Therefore, we can apply Theorem~\ref{theo:MT2} to conclude that
        \[
            \sum_{\Gamma \in \Sink^*(\mathfrak{u})} \, \prod_{x \in V(\Gamma)} \alpha_r(\updelta(x)) \,\leq\, \frac{\beta_r(\mathfrak{u})}{1 - \beta_r(\mathfrak{u})} \,=\, \frac{\eta}{1 - \eta}.
        \]
        Putting everything together,
        \begin{align*}
            \sum_{r < R} \, \sum_{\Gamma \in \mathbb{G}_r} p^{(1-\epsilon)(|V(\Gamma)| - 1)} 
            \,&\leq\, R \cdot \frac{1}{\xi\, (1+\eta)^{N}} \cdot \frac{\eta}{1 - \eta} \,=\, (1+\eta)^{-N} \,F(d, \eta, R),
        \end{align*}
        where $F(d, \eta, R)$ depends only on $d$, $\eta$, and $R$, and not on $N$, as desired.
    \end{scproof}

    \section{Solving Borel CSPs}\label{sec:Borel}

    \subsection{Preliminaries}\label{subsec:Borel_defns}

    In \S\ref{sec:Borel} we assume some basic familiarity with descriptive set theory. As mentioned in the introduction, the reader is directed to \cite{KechrisDST,AnushDST} for the necessary background. In particular, we shall make significant use of the material in \cite[Chap.~17 and 18]{KechrisDST}.
    
    In this subsection we give the formal definition of a Borel CSP. To begin with, for a standard Borel space $X$ and $n \in \N$, we let $[X]^n$ \ep{resp.~$\fins{X}$} be the set of all $n$-element \ep{resp.~finite} subsets of $X$. Both $[X]^n$ and $\fins{X}$ are naturally endowed with standard Borel structures. Namely, we let $(X)^n \subseteq X^n$ be the Borel set of all ordered $n$-tuples of pairwise distinct points and define an equivalence relation $\sim_n$ on $(X)^n$ by
    \[
        (x_0, \ldots, x_{n-1}) \,\sim_n\, (y_0, \ldots, y_{n-1}) \quad \Longleftrightarrow \quad \set{x_0, \ldots, x_{n-1}} \,=\, \set{y_0, \ldots, y_{n-1}}.
    \]
    By \cite[Ex.~6.1 and Prop.~6.3]{KechrisMiller}, the quotient space $(X)^n/{\sim_n} \cong [X]^n$ is standard Borel. Hence, the disjoint union $\fins{X} = [X]^0 \cup [X]^1 \cup [X]^2 \cup \ldots$ is standard Borel as well \cite[Prop.~1.4]{AnushDST}. 

    Next, given standard Borel spaces $X$ and $Y$, we let $\finf{X}{Y}$ be the set of all partial maps $\phi \colon X \pto Y$ with finite domains. Identifying each such map $\phi$ with its graph, i.e., with the set \[\set{(x,\phi(x)) \,:\, x \in \dom(\phi)} \,\subseteq\, X \times Y,\] makes $\finf{X}{Y}$ a Borel subset of $\fins{X \times Y}$ and hence a standard Borel space in its own right.
    
    Now we can formally define Borel CSPs:

    \begin{defn}[Borel CSPs]
        A CSP $\Pi = (V, \Lambda, \Const, \dom, \Bad)$ is \emphd{Borel} if:
        \begin{itemize}
            \item $V$, $\Lambda$, and $\Const$ are standard Borel spaces,
            \item the function $\dom \colon \Const \to \fins{V}$ is Borel, 
            \item the set $\set{(\const, \phi) \in \Const \times \finf{V}{\Lambda} \,:\, \phi \in \Bad(\const)}$ is Borel.
        \end{itemize}
    \end{defn}

    A graph $G$ is \emphd{Borel} if $V(G)$ is a standard Borel space and $E(G)$ is a Borel subset of $[V(G)]^2$. It is routine to check that if $\Pi$ is a Borel CSP, then its dependency graph $D_\Pi$ is Borel.




    For the remainder of \S\ref{sec:Borel}, we let $\Pi = (V, \Lambda, \Const, \dom, \Bad)$ be a Borel CSP and fix a Borel probability measure $\P$ on $\Lambda$. As in \S\ref{sec:MT}, we write $D \defeq D_\Pi$ and let $N[\const]$ be the closed neighborhood of a constraint $\const$ in $D$. Given a set $\mathfrak{U} \subseteq \Const$, we let $\dom(\mathfrak{U}) \defeq \bigcup_{\const \in \mathfrak{U}} \dom(\const)$ and $N[\mathfrak{U}] \defeq \bigcup_{\const \in \mathfrak{U}} N[\const]$. Note that if $D$ is locally finite and $\mathfrak{U} \subseteq \Const$ is Borel, then the sets $\dom(\mathfrak{U})$ and $N[\mathfrak{U}]$ are also Borel by the Luzin--Novikov theorem \cite[Thm.~18.10]{KechrisDST}.

    \subsection{Proof of Lemma~\ref{lemma:double_exp}}\label{subsec:proof_de}

    We start with a simple observation, which is a consequence of a standard ``large section'' uniformization result in descriptive set theory \cite[\S18.B]{KechrisDST}:

    \begin{Lemma}\label{lemma:edgeless}
        If $D$ has no edges and $\P[\Bad(\const)] < 1$ for all $\const \in \Const$, then $\Pi$ has a Borel solution.
    \end{Lemma}
    \begin{scproof}
        Since $\P[\Lambda^{\dom(\const)} \setminus \Bad(\const)] > 0$ for all $\const \in \Const$, by \cite[Corl.~18.7]{KechrisDST}, there is a Borel function
        \[
            F \colon \Const \to \finf{V}{\Lambda}
        \]
        such that $F(\const) \in \Lambda^{\dom(\const)} \setminus \Bad(\const)$ for every constraint $\const \in \Const$. Since $D$ has no edges \ep{i.e., the constraints in $\Const$ have disjoint domains}, we can pick an arbitrary value $\lambda \in \Lambda$ and define, for each $v \in V$,
        \[
            f(v) \,\defeq\, \begin{cases}
                \big(F(\const)\big)(v) &\text{if } v \in \dom(\const) \text{ for some \ep{unique} } \const \in \Const,\\
                \lambda &\text{if } v \notin \dom(\Const).
            \end{cases}
        \]
        Then $f \colon V \to \Lambda$ is a Borel solution to $\Pi$, as desired. 
    \end{scproof}
    
    As mentioned in the introduction, the proof technique we employ to establish Lemma~\ref{lemma:double_exp} is the method of conditional probabilities, which is essentially based on a greedy algorithm. We will construct a Borel solution to the given CSP $\Pi$  in a finite number of stages. At the start of each stage, some of the variables will already have their labels  determined, and our goal will be to extend this partial labeling to a larger set of variables, yielding a solution to $\Pi$ in the end. We will achieve this by repeatedly applying Lemma~\ref{lemma:edgeless} to certain auxiliary CSPs.
    
    To help with the analysis of this inductive construction, we introduce the following notation. Given $U \subseteq V$ and $f \colon U \to \Lambda$, the \emphd{quotient CSP} $\Pi/f$ encodes the problem of extending $f$ to a full solution to $\Pi$. That is, we define
    \[
        \Pi/f \,\defeq\, \big(V \setminus U, \, \Lambda, \, \Const, \, \dom/f, \, \Bad/f \big),
    \]
    where, for $\const \in \Const$, $(\dom/f)(\const) \defeq \dom(\const) \setminus U$ and $(\Bad/f)(\const)$ is the set of all labelings $\phi   \in \Lambda^{\dom(\const) \setminus U}$ such that extending $f$ by $\phi$ would violate $\const$. In other words, for all $\phi \colon (\dom(\const) \setminus U) \to \Lambda$, we have
    \[
        \phi \,\in\, (\Bad/f)(\const) \quad \Longleftrightarrow \quad \rest{f}{\dom(\const) \cap U} \,\cup\, \phi \,\in\, \Bad(\const).
    \]
    By definition, $f' \colon (V \setminus U) \to \Lambda$ is a solution to $\Pi/f$ if and only if $f \cup f'$ is a solution to $\Pi$. Let us now make a few observations:
    
        \begin{enumerate}[wide, label=\ep{\normalfont\roman*}]
            \item It is possible that $\dom(\const) \subseteq U$ for some $\const \in \Const$. In this case, $(\dom/f)(\const) = \0$ and $(\Bad/f)(\const)$ is equal to $\0$ if $f$ satisfies $\const$  and $\set{\0}$ if $f$ violates $\const$ \ep{in the CSP $\Pi$}.


        \item Since $(\dom/f)(\const) \subseteq \dom(\const)$ for all $\const \in \Const$, the dependency graph of $\Pi/f$ is a subgraph of $D$. In particular, its maximum degree is at most that of $\Pi$. 

        \item If $U$ is a Borel subset of $V$ and $f$ is a Borel function, then $\Pi/f$ is a Borel CSP.
        \end{enumerate}

    We can now describe the inductive step in the proof of Lemma~\ref{lemma:double_exp}. For a graph $G$ and $r \in \N$, let $G^r$ be the graph with 
    $V(G^r) \defeq V(G)$ and $E(G^r) \defeq \set{\set{x,y} \in [V(G)]^2 \,:\, y \in B_G(x,r)}$. 

    \begin{Lemma}\label{lemma:induction}
        Suppose the maximum degree of $D$ is at most $d \in \N$. If $\I \subseteq \Const$ is a Borel $D^2$-independent set, then there is a Borel function $f \colon \dom(\I) \to \Lambda$ such that for all $\mathfrak{a} \in \Const$,
        \begin{equation}\label{eq:step}
            \P[(\Bad/f)(\mathfrak{a})] \,\leq\, \begin{cases}
                (d+1) \,\P[\Bad(\mathfrak{a})] &\text{if } \mathfrak{a} \in N[\I],\\
                \P[\Bad(\mathfrak{a})] &\text{if } \mathfrak{a} \notin N[\I].
            \end{cases}
        \end{equation}
    \end{Lemma}
    \begin{scproof}
    Form an auxiliary CSP $\Pi'$ as follows:
        \[
            \Pi' \,\defeq\, \big(\dom(\I), \, \Lambda, \, \I, \, \rest{\dom}{\I}, \, \mathcal{A}\big),
        \]
        where for each $\const \in \I$, $\mathcal{A}(\const)$ is the set of all mappings $\phi \colon \dom(\const) \to \Lambda$ such that 
        \[
            \P[(\Bad/\phi)(\mathfrak{a})] \,>\, (d+1) \,\P[\Bad(\mathfrak{a})] \quad \text{for some } \mathfrak{a} \in N[\const].
        \]
        Note that the set
        \[
            Q \,\defeq\, \big\{(\const, \phi, \mathfrak{a}) \in \I \times \finf{V}{\Lambda} \times \Const \,:\, \phi \in \Lambda^{\dom(\const)} \text{ and } \P[(\Bad/\phi)(\mathfrak{a})] > (d+1) \,\P[\Bad(\mathfrak{a})] \big\}
        \]
        is Borel by \cite[Thm.~17.25]{KechrisDST}. Therefore, the set
        \begin{align*}
            \big\{(\const, \phi) \in \I \times \finf{V}{\Lambda} \,:\, \phi \in \mathcal{A}(\const)\big\} \,=\, \big\{(\const, \phi) \,:\, \exists \mathfrak{a} \in N[\const] \, \big((\const, \phi, \mathfrak{a}) \in Q\big)\big\}
        \end{align*}
        is also Borel by the Luzin--Novikov theorem \cite[Thm.~18.10]{KechrisDST}. It follows that $\Pi'$ is a Borel CSP.

        Since $\I$ is $D^2$-independent, the constraints in $\I$ have disjoint domains, so the dependency graph of $\Pi'$ has no edges. Next we argue that 
        $\P[\mathcal{A}(\const)] < 1$ for all $\const \in \I$. 
%
            For $\mathfrak{a} \in N[\const]$ and $\phi \in \Lambda^{\dom(\const)}$, let 
            \[
                X_\mathfrak{a}(\phi) \,\defeq\, \P[(\Bad/\phi)(\mathfrak{a})].
            \]
            By Fubini's theorem, the expectation of the random variable $X_\mathfrak{a}(\phi)$ with respect to a random choice of $\phi \sim \P^{\dom(\const)}$ is $\E[X_\mathfrak{a}(\phi)] = \P[\Bad(\mathfrak{a})]$. 
            If $\P[\Bad(\mathfrak{a})] = 0$, then $X_\mathfrak{a}(\phi) = 0$ almost surely. On the other hand, if $\P[\Bad(\mathfrak{a})] > 0$, then, by Markov's inequality,
            \begin{align*}
                \P\big[X_\mathfrak{a}(\phi) > (d+1) \,\P[\Bad(\mathfrak{a})]\big] \,<\, \frac{\E[X_\mathfrak{a}(\phi)]}{(d+1)\,\P[\Bad(\mathfrak{a})]} \, =\, \frac{1}{d+1}.
            \end{align*}
            Therefore, regardless of the value of $\P[\Bad(\mathfrak{a})]$, the union bound gives
            \[
                \P[\mathcal{A}(\const)] \,\leq\, \sum_{\mathfrak{a} \in N[\const]} \P\big[X_\mathfrak{a}(\phi) > (d+1) \,\P[\Bad(\mathfrak{a})]\big] \,<\, \frac{1}{d+1} \cdot |N[\const]| \,\leq\, 1. 
            \]

         It follows that $\Pi'$ has a Borel solution $f \colon \dom(\I) \to \Lambda$ by Lemma~\ref{lemma:edgeless}. We claim that $f$ is as desired. Indeed, take any $\mathfrak{a} \in \Const$. We verify \eqref{eq:step} for $\mathfrak{a}$. If $\mathfrak{a} \notin N[\I]$, then $(\Bad/f)(\mathfrak{a}) = \Bad(\mathfrak{a})$ and we are done. On the other hand, if $\mathfrak{a} \in N[\I]$, then there is a constraint $\const \in \I$ with $\mathfrak{a} \in N[\const]$, and such $\const$ is unique because $\I$ is $D^2$-independent. Let $\phi \defeq \rest{f}{\dom(\const)}$. Then $(\Bad/f)(\mathfrak{a}) = (\Bad/\phi)(\mathfrak{a})$. Since $f$ is a solution to $\Pi'$, we have $\phi \notin \mathcal{A}(\const)$, which yields
        $
            \P[(\Bad/\phi)(\mathfrak{a})] \leq (d+1) \, \P[\Bad(\mathfrak{a})]
        $, as claimed.
    \end{scproof}

    \begin{scproof}[ of Lemma~\ref{lemma:double_exp}]
        Suppose $p \in [0,1)$ and $d \in \N$ are such that $\P[\Bad(\const)] \leq p$ for all $\const \in \Const$, the maximum degree of $D$ is at most $d$, and $p\, (d+1)^{d+1} < 1$. We seek a Borel solution to $\Pi$.

        Note that the maximum degree of $D^2$ is at most $d^2$.  Hence, by \cite[Prop.~4.6]{kechris1999borel}, there is a partition $\Const = \I_1 \cup \ldots \cup \I_n$ of $\Const$ into $n \leq d^2 + 1$ Borel $D^2$-independent sets. Let \[U_i \,\defeq\, \dom(\I_i) \setminus (\dom(\I_1) \cup \ldots \cup \dom(\I_{i-1})).\] Applying Lemma~\ref{lemma:induction} iteratively, we obtain a sequence of Borel labelings $f_i \colon U_i \to \Lambda$ such that for all $\mathfrak{a} \in \Const$ and $1 \leq i \leq n$,
        \[
            \P[(\Bad/(f_1 \cup \ldots \cup f_i))(\mathfrak{a})] \,\leq\, \begin{cases}
                (d+1) \,\P[(\Bad/(f_1 \cup \ldots \cup f_{i-1}))(\mathfrak{a})] &\text{if } \mathfrak{a} \in N[\I_i],\\
                \P[(\Bad/(f_1 \cup \ldots \cup f_{i-1}))(\mathfrak{a})] &\text{if } \mathfrak{a} \notin N[\I_i].
            \end{cases}
        \]
        Fix any $\lambda \in \Lambda$ and define $f \colon V \to \Lambda$ by
        \[
            f(v) \,\defeq\, \begin{cases}
                f_i(v) &\text{if } v \in U_i \text{ for some $1 \leq i \leq n$},\\
                \lambda &\text{if } v\notin \dom(\Const).
            \end{cases}
        \]
        Since for each constraint $\mathfrak{a} \in \Const$, there are $|N[\mathfrak{a}]|$ indices $1 \leq i \leq n$ with $\mathfrak{a} \in N[\I_i]$, we have
        \[
            \P[(\Bad/f)(\mathfrak{a})] \,=\, \P[(\Bad/(f_1 \cup \ldots \cup f_n))(\mathfrak{a})] \,\leq\, (d+1)^{|N[\mathfrak{a}]|} \,\P[\Bad(\mathfrak{a})] \,\leq\, (d+1)^{d+1} \,p \,<\, 1.
        \]
        But $f$ is defined on all of $V$, so $(\dom/f)(\mathfrak{a}) = \0$, i.e.,  $(\Bad/f)(\mathfrak{a})$ is either $\0$ or $\set{\0}$. The latter case is impossible as $\P[\set{\0}] = 1$, hence $(\Bad/f)(\mathfrak{a}) = \0$. This means that $f$ satisfies $\mathfrak{a}$. Therefore, $f$ is a Borel solution to $\Pi$, and the proof is complete.
    \end{scproof}

    \subsection{Proof of Theorem~\ref{theo:main}}\label{subsec:final}

    Recall the set-up of Theorem~\ref{theo:main}: We have 
    parameters $p \in [0,1)$, $d \in \N$, and $s > 1$ such that:
        \begin{enumerate}[label=\ep{\normalfont\roman*}]
            \item $\P[\Bad(\const)] \leq p$ for all $\const \in \Const$ and the maximum degree of $D$ is at most $d$,
            \item $\egr(D) < s$, and
            \item $p\, (\e \,(d+1))^{s} < 1$.
        \end{enumerate}
    Our goal is to find a Borel solution to $\Pi$.

    We start by observing that the Moser--Tardos Algorithm can be used in the Borel setting:

    \begin{claim}\label{claim:good_Borel}
        If $\tbl \colon V \to \Lambda^\N$ is a Borel good table, then $\Pi$ has a Borel solution $f \in \MMTA(\Pi, \tbl)$.
    \end{claim}
    \begin{claimproof}
        It is enough to argue that there is a Borel function $f \in \MMTA(\Pi, \tbl)$, since such $f$ would be defined on all of $V$ because $\tbl$ is good, and hence be a solution to $\Pi$ by Lemma~\ref{lemma:limit}. To this end, we observe that if $\Pi$ and $\tbl$ are Borel, then the Maximal Moser--Tardos Algorithm can be carried out on input $(\Pi,\tbl)$ so that all of $\level_n$, $f_n$, $\Const_n$, $\I_n$, $\level$, $f$ are Borel. Indeed, $\level_0$ is Borel because so is $\tbl$. If $\level_n$ is Borel, then so is $f_n$, and $\Const_n$ is Borel as well because $\Pi$ is Borel. We can then choose a Borel maximal $D$-independent subset $\I_n \subseteq \Const_n$ using \cite[Prop.~4.2]{kechris1999borel}. The set $\dom(\I_n)$ is also Borel, and thus 
        $\level_{n+1}$ is Borel as well. Finally, if $\level_n$ is Borel for all $n \in \N$, then $\level$ and $f$ are also Borel. 
    \end{claimproof} 

    Next we show that the CSP $\mathsf{LG}(R, N, \epsilon)$ defined in \S\ref{subsec:good_definitions} (specifically in Lemma~\ref{lemma:local_CSP}) is Borel: 

    \begin{claim}\label{claim:Borel}
        For each $R$, $N \in \N$ and $\epsilon \in (0,1)$, the CSP $\mathsf{LG}(R, N, \epsilon)$ is Borel.
    \end{claim}
    \begin{claimproof}
        The function $\dom_R \colon \Const \to \fins{V}$ is Borel by the Luzin--Novikov theorem \cite[Thm.~18.10]{KechrisDST}, see \cite[Corl.~5.2]{pikhurko2020borel} for details. To see that the set $Q \defeq \set{(\const, \phi) \in \Const \times \finf{V}{\Lambda^\N} \,:\, \phi \in \LBad_{R,N,\epsilon}(\const)}$ is Borel, we use the characterization of $(\const, R, N, \epsilon)$-locally good tables from Lemma~\ref{lemma:restricted_Folner}, which shows that, as there are only countably many $(\const, R)$-bounded finite $\MT$-sequences for each $\const \in \Const$ and $R \in \N$, $Q$ is also Borel by the Luzin--Novikov theorem \cite[Thm.~18.10]{KechrisDST}.
    \end{claimproof}

    Take $\epsilon \in (0,1)$ such that $\egr(D) < (1-\epsilon)^{-1} < s$. By the definition of the exponential growth rate, there is $R \in \N$ such that $\gamma_D(R) < (1 - \epsilon)^{-R}$. On the other hand, the inequality $(1-\epsilon)^{-1} < s$ is equivalent to $\epsilon + 1/s < 1$, and we can take $\eta \in (0,1)$ so small that
    \[
        \frac{1-\eta}{1+\eta} \,>\, p^{1 - \epsilon - \frac{1}{s}}.
    \]
    Finally, we take $N \in \N$ so large that
    \[
        (1+\eta)^{N} \, > \, F(d, \eta, R) \, \gamma_D(2R+1)^{\gamma_D(2R+1)},
    \]
    where $F(d, \eta, R)$ is the function from Lemma~\ref{lemma:LG_prob}. Claim~\ref{claim:Borel} and Lemmas~\ref{lemma:LG_degree} and \ref{lemma:LG_prob} imply that the CSP $\mathsf{LG}(R, N, \epsilon)$ satisfies the assumptions of Lemma~\ref{lemma:double_exp}. Therefore, $\mathsf{LG}(R, N, \epsilon)$ has a Borel solution $\tbl \colon V \to \Lambda^\N$. In other words, $\tbl$ is a Borel $(R, N,\epsilon)$-locally good table. By Lemma~\ref{lemma:good}, $\tbl$ is good, so, by Claim~\ref{claim:good_Borel}, there exists a Borel solution to $\Pi$, and we are done.